\documentclass[twocolumn]{article}
\usepackage{epsfig}
\usepackage{graphicx}
\usepackage{color}
\usepackage[acmtitlespace,nocopyright]{acmneww}
\usepackage{times}
\usepackage{amsmath}
\usepackage{amssymb}
\usepackage{xspace}

\newcommand {\mm}[1] {\ifmmode{#1}\else{\mbox{\(#1\)}}\fi}

\newcommand {\scalprod}[2] {{\langle #1 , #2 \rangle}}

\newcommand{\proof}{\noindent{\sc Proof.~}}

\newcommand{\eop}{\hfill\usebox{\smallProofsym}\bigskip}  %
\newsavebox{\smallProofsym}                            
\savebox{\smallProofsym}                               %
{
\begin{picture}(6,6)                                   %
\put(0,0){\framebox(6,6){}}                            %
\put(0,2){\framebox(4,4){}}                            %
\end{picture}                                          %
}                                                      %
\makeatletter
\long\def\@makecaption#1#2{%
  \vskip\abovecaptionskip
  \sbox\@tempboxa{\small #1: #2}%
  \ifdim \wd\@tempboxa >\hsize
    \small #1: #2\par
  \else
    \global \@minipagefalse
    \hb@xt@\hsize{\hfil\box\@tempboxa\hfil}%
  \fi
  \vskip\belowcaptionskip}
\makeatother

\newcommand{\Bgroup}        {\mm{\sf B}}
\newcommand{\Cgroup}        {\mm{\sf C}}
\newcommand{\Dgroup}        {\mm{\sf D}}

\newcommand{\Fgroup}        {\mm{\sf F}}
\newcommand{\Hgroup}        {\mm{\sf H}}
\newcommand{\rHgroup}       {\mm{\sf \tilde{H}}}
\newcommand{\Mgroup}        {\mm{\sf M}}
\newcommand{\Sgroup}        {\mm{\sf S}}

\newcommand{\rUgroup}       {\mm{\sf \tilde{U}}}

\newcommand{\rVgroup}       {\mm{\sf \tilde{V}}}
\newcommand{\Xgroup}        {\mm{\sf X}}
\newcommand{\rXgroup}       {\mm{\sf \tilde{X}}}
\newcommand{\Ygroup}        {\mm{\sf Y}}
\newcommand{\Zgroup}        {\mm{\sf Z}}

\newcommand{\gmap}          {\mm{\sf g}}
\newcommand{\hmap}          {\mm{\sf h}}
\newcommand{\isomap}        {\mm{\sf iso}}
\newcommand{\mmap}          {\mm{\sf m}}
\newcommand{\ik}            {\mm{\sf ik}}
\newcommand{\ip}            {\mm{\sf ip}}
\newcommand{\ck}            {\mm{\sf ck}}
\newcommand{\cp}            {\mm{\sf cp}}

\newcommand{\Betti}         {\mm{\beta}}
\newcommand{\rBetti}        {\mm{\tilde{\beta}}}
\newcommand{\Ddgm}          {\mm{\sf Dgm}}

\newcommand{\rDdgm}         {\mm{\sf \tilde{D}gm}}

\newcommand{\Odgm}          {\mm{\sf Ord}}
\newcommand{\Hdgm}          {\mm{\sf Hor}}
\newcommand{\Vdgm}          {\mm{\sf Vcl}}
\newcommand{\Rdgm}          {\mm{\sf Rel}}
\newcommand{\Lwing}         {\mm{\cal L}}
\newcommand{\Rwing}         {\mm{\cal R}}
\newcommand{\rLwing}        {\mm{\cal \tilde{L}}}
\newcommand{\rRwing}        {\mm{\cal \tilde{R}}}

\newcommand{\kernel}[1]     {\mm{\sf ker\,}{#1}}
\newcommand{\image}[1]      {\mm{\sf im\,}{#1}}
\newcommand{\pimage}[1]     {\mm{\sf pim\,}{#1}}
\newcommand{\cokernel}[1]   {\mm{\sf cok\,}{#1}}

\newcommand{\Aspace}        {\mm{{\mathbb A}}}
\newcommand{\Bspace}        {\mm{{\mathbb B}}}
\newcommand{\Hspace}        {\mm{{\mathbb H}}}
\newcommand{\Mspace}        {\mm{{\mathbb M}}}
\newcommand{\Rspace}        {\mm{{\mathbb R}}}
\newcommand{\Sspace}        {\mm{{\mathbb S}}}
\newcommand{\Uspace}        {\mm{{\mathbb U}}}
\newcommand{\Vspace}        {\mm{{\mathbb V}}}
\newcommand{\Xspace}        {\mm{{\mathbb X}}}
\newcommand{\Yspace}        {\mm{{\mathbb Y}}}
\newcommand{\Zspace}        {\mm{{\mathbb Z}}}

\newcommand{\rank}[1]       {\mm{\rm rank\,}{#1}}

\newcommand{\capsp}         {{\; \cap \;}}
\newcommand{\cupsp}         {{\; \cup \;}}
\newcommand{\from}          {\mm{\leftarrow}}

\newtheorem{result}{}

\newcommand{\ignore}[1]{}

\title{Alexander Duality for Functions: ~\\
       the Persistent Behavior of Land and Water and Shore
       \thanks{This research is partially supported by
               the National Science Foundation (NSF) under grant DBI-0820624,
               and the European Science Foundation under
               the Research Networking Programme.}
       }

\author{Herbert Edelsbrunner\thanks{IST Austria (Institute of Science and
            Technology Austria), Kloster\-neu\-burg, Austria,
            Departments of Computer Science and of Mathematics,
            Duke University, Durham, North Carolina,
            and Geomagic, Research Triangle Park, North Carolina.}
        and Michael Kerber\thanks{IST Austria (Institute of Science and
            Technology Austria), Kloster\-neu\-burg, Austria.}}

\begin{document}
\maketitle

\begin{abstract}
  This note contributes to the point calculus of persistent homology
  by extending Alexander duality to real-valued functions.
  Given a perfect Morse function $f: \Sspace^{n+1} \to [0,1]$
  and a decomposition $\Sspace^{n+1} = \Uspace \cupsp \Vspace$
  such that $\Mspace = \Uspace \capsp \Vspace$ is an $n$-manifold,
  we prove elementary relationships between the persistence diagrams
  of $f$ restricted to $\Uspace$, to $\Vspace$, and to $\Mspace$.
\end{abstract}

\vspace{0.1in}
{\small
 \noindent{\bf Keywords.}
   Algebraic topology, homology, Alexander duality,
   Mayer-Vietoris sequences, persistent homology, point calculus.}

\section{Introduction}
\label{sec1}

Persistent homology is a recent extension of the classical theory
of homology; see e.g.\ \cite{EdHa10}.
Given a real-valued function on a topological space,
it measures the importance of a homology class by monitoring
when the class appears and when it disappears in the increasing
sequence of sublevel sets.
A technical requirement is that the function be \emph{tame},
which means it has only finitely many homological critical values,
and each sublevel set has finite rank homology groups.
Pairing up the births and deaths, and drawing each pair of values
as a \emph{dot} (a point in the plane), we get a multiset which we
refer to as the \emph{persistence diagram} of the function.
It is a combinatorial summary of the homological information
contained in the sequence of sublevel sets.
If we substitute reduced for standard homology groups,
we get a slightly modified \emph{reduced persistence diagram}
of the function.

As between homology groups, we can observe relationships between
persistence diagrams.
A prime example is Poincar\'{e} duality, which says that
the $p$-th and the $(n-p)$-th homology groups of an $n$-manifold
are isomorphic.
More precisely, this is true if homology is defined for field
coefficients, which is what we assume throughout this paper.
The extension to functions says that the diagram
is symmetric with respect to reflection across the vertical axis;
see \cite{CEH09}.
Here, we change the homological dimension of a dot from
$p$ to $n-p$ whenever we reflect it across the axis.
This paper contributes new relationships by extending
Alexander duality from spaces to functions.
To state our results, we assume a perfect Morse function,
$f: \Sspace^{n+1} \to [0,1]$, which for the sphere
has no critical points other than a minimum and a maximum,
and a decomposition of the $(n+1)$-dimensional sphere into two subsets,
$\Sspace^{n+1} = \Uspace \cupsp \Vspace$,
whose intersection is an $n$-manifold, $\Mspace = \Uspace \capsp \Vspace$.
Our first result says that the reduced persistence diagrams
of $f$ restricted to $\Uspace$ and to $\Vspace$ are reflections of each other.
We call this the Land and Water Theorem.
Our second result relates land with shore.
Ignoring some modifications, it says that the persistence diagram
of $f$ restricted to $\Mspace$ is the disjoint union of the diagram
of $f$ restricted to $\Uspace$ and of its reflection.
The modifications become unnecessary if we assume that the
minimum and maximum of $f$ both belong to a common component of $\Vspace$.
We call this the Euclidean Shore Theorem.

In the example that justifies the title of this paper, and the names
of our theorems,
we let $\Uspace$ be the planet Earth, not including the water and the air.
To a coarse approximation, $\Uspace$ is homeomorphic to a $3$-ball,
sitting inside the Universe, which we model as a $3$-sphere, $\Sspace$.
The function we consider is the negative gravitational potential of the Earth,
which is defined on the entire Universe.
The sea is then a sublevel set of this function restricted to $\Vspace$,
which is the closure of $\Sspace - \Uspace$.
With these definitions, our results relate the persistence diagram
of the gravitational potential restricted to the Earth with the
shape of the sea as its water level rises.
Also, the Euclidean Shore Theorem applies, unmodified,
expanding the relationship to include the sea floor, 
which is swept out by the shoreline as the water level rises.

Besides developing the mathematical theory of persistent homology,
there are pragmatic reasons for our interest in the
extension of Alexander duality to functions.
Persistence has fast algorithms, so that the bulk of the work is
often in the construction of the space and the function for which
we compute persistence.
A point in case is the analysis of the biological process
of cell segregation started in \cite{EHKKS11}.
Modeling the process as a subset of space-time, the function
of interest is time which, after compactifying space-time
to $\Sspace^4$, has no critical points other than a minimum and a maximum.
The subset $\Uspace$ of $\Sspace^4$ is a union of cells
times time, whose boundary is a $3$-manifold.
We can represent $\Uspace$ by a $1$-parameter family of alpha complexes,
whose disjoint union has the same homotopy type;
see e.g.\ \cite[Chapter III]{EdHa10}.
However, the boundary of that disjoint union is not necessarily
a $3$-manifold.
Using our Euclidean Shore Theorem, we can compute the persistence diagram
of the function on the $3$-manifold without ever constructing the $3$-manifold.

\paragraph{Outline.}
Sections \ref{sec2} and \ref{sec3} introduce the necessary background
on homology and persistent homology.
Sections \ref{sec4} and \ref{sec5} present our two results.
Section \ref{sec6} concludes the paper.

\section{Homology}
\label{sec2}

Starting with a brief introduction of classical homology groups,
we present the relevant background on Alexander duality and
Mayer-Vietoris sequences.
More comprehensive discussions of these topics can be found
in textbooks of algebraic topology, such as \cite{Hat02,Mun84}.

\paragraph{Background.}
The $p$-dimensional homology of a topological space, $\Xspace$,
is a mathematical language to define, count, and reason about
the $p$-dimensional connectivity of $\Xspace$.
There are different but essentially equivalent theories
depending on the choices one makes in the representation of
the space, the selection of cycles, and the meaning of addition.
For our purpose, the most elementary of these theories
will suffice:
a simplicial complex, $K$, that triangulates $\Xspace$,
formal sums of $p$-simplices with zero boundary as $p$-cycles,
and adding with coefficients in a field, $\Fgroup$.
Most algorithms on homology assume this model,
in particular the ones developed within persistent homology,
which is defined only for field coefficients.
In this model, we call a formal sum of $p$-simplices a
\emph{$p$-chain}, a \emph{$p$-cycle} if its boundary is empty,
and a \emph{$p$-boundary} if it is the boundary of a $(p+1)$-chain.
The $p$-boundaries form a subgroup of the $p$-cycles,
which form a subgroup of the $p$-chains:
$\Bgroup_p \subseteq \Zgroup_p \subseteq \Cgroup_p$.
The \emph{$p$-th homology group} is the quotient of
the $p$-cycles over the $p$-boundaries:
$\Hgroup_p = \Zgroup_p / \Bgroup_p$.
Its elements are sets of $p$-cycles that differ from each
other by $p$-boundaries.
More fully, we denote the $p$-th homology group
by $\Hgroup_p (K, \Fgroup)$,
or by $\Hgroup_p (\Xspace, \Fgroup)$ to emphasize that
the group is independent of the simplicial complex
we choose to triangulate the space.
However, we will fix an arbitrary field $\Fgroup$ and drop it from the notation.
For field coefficients, $\Hgroup_p (\Xspace)$ is necessarily
a vector space, which is fully described by its rank,
$\Betti_p (\Xspace) = \rank{\Hgroup_p (\Xspace)}$
such that $\Hgroup_p (\Xspace) \simeq \Fgroup^{\Betti_p (\Xspace)}$.
This rank is called the \emph{$p$-th Betti number} of $\Xspace$.
Finally, we will often drop the dimension from the notation
by introducing the direct sum,
$\Hgroup (\Xspace) = \bigoplus_{p \in \Zspace} \Hgroup_p (\Xspace)$.

Besides standard homology, we will frequently use
\emph{reduced homology groups}, $\rHgroup_p (\Xspace)$,
which are isomorphic to the non-reduced groups except possibly
for dimensions $p = 0, -1$.
To explain the difference, we note that $\Betti_0 (\Xspace)$
counts the components of $\Xspace$.
In contrast, $\rBetti_0 (\Xspace) = \rank{\rHgroup_0 (\Xspace)}$
counts the gaps between components or, equivalently,
the edges that are needed to merge the components into one.
Hence, $\rBetti_0 (\Xspace) = \Betti_0 (\Xspace) - 1$,
except when $\Xspace$ is empty, in which case
$\rBetti_0 (\Xspace) = \Betti_0 (\Xspace) = 0$.
To distinguish this case from a single component,
we have $\rBetti_{-1} (\Xspace)$ equal to $1$
if $\Xspace = \emptyset$, and equal to $0$, otherwise.
Furthermore, we use relative homology, which is defined for
pairs of spaces, $\Yspace \subseteq \Xspace$.
Taking a pair relaxes the requirement of a chain to be called a cycle,
namely whenever its boundary is contained in $\Yspace$,
which includes the case when the boundary is empty.
We write $\Hgroup_p (\Xspace, \Yspace)$
for the \emph{$p$-th relative homology group} of the pair,
and $\Betti_p (\Xspace, \Yspace) = \rank{\Hgroup_p (\Xspace, \Yspace)}$
for the \emph{$p$-th relative Betti number}.
As before, we will suppress the dimension from the notation
by introducing
$\Hgroup (\Xspace, \Yspace) = \bigoplus_{p \in \Zspace}
                                \Hgroup_p (\Xspace, \Yspace)$.

As examples, consider the $(n+1)$-dimensional sphere,
$\Sspace = \Sspace^{n+1}$,
and a closed hemisphere, $\Hspace \subseteq \Sspace$,
which is a ball of dimension $n+1$.
In standard homology, we have
$\Betti_0 (\Sspace) = \Betti_0 (\Hspace) = \Betti_{n+1} (\Sspace) = 1$
while all other Betti numbers are zero.
In reduced homology, we have
$\rBetti_{n+1} (\Sspace) = 1$
while all other reduced Betti numbers are zero.
In particular, $\rBetti_p (\Hspace) = 0$ for all $p$.
In relative homology, we have
$\Betti_{n+1} (\Sspace, \Hspace) = 1$
while all other relative Betti numbers are zero.
In particular, $\Betti_0 (\Sspace, \Hspace) = 0$,
which may be confusing at first but makes sense because
every point on the sphere can be connected by a path to a
point in the hemisphere and is thus a $0$-boundary.

\paragraph{Alexander duality.}
Recall that a \emph{perfect Morse function} is one whose number
of critical points equals the sum of Betti numbers of the space.
For a sphere, this number is $2$:  a minimum and a maximum.
Throughout the remainder of this section,
we assume a perfect Morse function $f: \Sspace \to [0,1]$,
whose values at the minimum and the maximum are $0$ and $1$.
We also assume \emph{complementary subsets} whose union is
the sphere, $\Uspace \cupsp \Vspace = \Sspace$, and whose intersection
is an $n$-manifold, $\Uspace \capsp \Vspace = \Mspace$.
We assume that the restriction of $f$ to $\Mspace$ is tame,
which implies that its restrictions to $\Uspace$ and $\Vspace$ are also tame.
For each $t \in \Rspace$, we write $\Sspace_t = f^{-1} [0,t]$
for the \emph{sublevel set} of $f$,
and $\Uspace_t = \Uspace \capsp \Sspace_t$,
$\Vspace_t = \Vspace \capsp \Sspace_t$,
$\Mspace_t = \Mspace \capsp \Sspace_t$
for the sublevel sets of the restrictions of $f$
to $\Uspace$, $\Vspace$, $\Mspace$.
Similarly, we write $\Sspace^t = f^{-1} [t,1]$ for the
\emph{superlevel set}, and $\Uspace^t = \Uspace \capsp \Sspace^t$,
$\Vspace^t = \Vspace \capsp \Sspace^t$,
$\Mspace^t = \Mspace \capsp \Sspace^t$
for the superlevel sets of the restrictions.

The basic version of Alexander duality is a statement about two
complementary subsets of the sphere;
see \cite[page 424]{Mun84}.
More specifically, it states that $\Hgroup_q (\Uspace)$
and $\Hgroup_p (\Vspace)$ are isomorphic, where $q = n - p$,
except for $p=0$ and $q = 0$ when the $0$-dimensional group
has an extra generator.
This implies
\begin{eqnarray}
  \Betti_0 (\Vspace)  &=&  \Betti_n (\Uspace) + 1 ,
  \label{eqn:ADSphere-0}                                   \\
  \Betti_p (\Vspace)  &=&  \Betti_q (\Uspace) ,
  \label{eqn:ADSphere}                                     \\
  \Betti_n (\Vspace)  &=&  \Betti_0 (\Uspace) - 1 ,
  \label{eqn:ADSphere-n}
\end{eqnarray}
for $1 \leq p \leq n-1$.
We will also need the version that deals with two complementary
subsets of the $(n+1)$-dimensional ball, $\Bspace^{n+1}$;
see \cite[page 426]{Mun84}.
Let $t\in(0,1)$ denote a regular value of $f|_{\Mspace}$.
Note first that $\Sspace_t$ is homeomorphic to $\Bspace^{n+1}$.
By excision, the homology groups of $\Uspace_t$ relative
to $\Uspace_t \capsp f^{-1} (t)$ are isomorphic to those
of $\Uspace$ relative to $\Uspace^t$.
Alexander duality states that
$\Hgroup_q (\Uspace, \Uspace^t)$ and $\Hgroup_p (\Vspace_t)$
are isomorphic, where $q+p = n$, as before,
except for $p = 0$, when $\Vspace_t$ has an extra component.
This implies
\begin{eqnarray}
  \Betti_0 (\Vspace_t)  &=&  \Betti_n (\Uspace, \Uspace^t) + 1,
  \label{eqn:ADBall-0}                                              \\
  \Betti_p (\Vspace_t)  &=&  \Betti_q (\Uspace, \Uspace^t) ,
  \label{eqn:ADBall}
\end{eqnarray}
for $1 \leq p \leq n$.

\paragraph{Mayer-Vietoris.}
We can connect the homology groups of $\Uspace$ and $\Vspace$
with those of $\Mspace$ and $\Sspace$ using the 
Mayer-Vietoris sequence of the decomposition.
This sequence is \emph{exact}, meaning the image of every map
equals the kernel of the next map.
Counting the trivial homology groups, the sequence is infinitely long,
with three terms per dimension:

\begin{eqnarray*}
  \ldots  \to  \Hgroup_{p+1} (\Sspace)
          \to  \Hgroup_p (\Mspace)              
          \to  \Hgroup_p (\Uspace) \oplus \Hgroup_p (\Vspace) \\
          \to  \Hgroup_p (\Sspace)
          \to  \ldots ,
\end{eqnarray*}

where the maps between homology groups of the same dimension
are induced by the inclusions.
The only non-trivial homology groups of $\Sspace$
are in dimensions $0$ and $n+1$,
with ranks $\Betti_0 (\Sspace) = \Betti_{n+1} (\Sspace) = 1$.
It follows that for $1 \leq p \leq n-1$, the groups defined by
$\Mspace$, $\Uspace$, $\Vspace$ are surrounded by trivial groups.
This implies that the groups of $\Mspace$ and the direct sums
of the groups of $\Uspace$ and $\Vspace$ are isomorphic.
For $p = 0$ and $p = n$, the non-trivial homology groups of 
$\Sspace$ prevent this isomorphism, and we get
\begin{eqnarray}
  \Betti_0 (\Mspace)  &=&  \Betti_0 (\Uspace) + \Betti_0 (\Vspace) - 1,
    \label{eqn:MVSphere-0}                                                \\
  \Betti_p (\Mspace)  &=&  \Betti_p (\Uspace) + \Betti_p (\Vspace),
    \label{eqn:MVSphere}                                                 \\
  \Betti_n (\Mspace)  &=&  \Betti_n (\Uspace) + \Betti_n (\Vspace) + 1,
    \label{eqn:MVSphere-n}
\end{eqnarray}
for $1 \leq p \leq n-1$. We also have a Mayer-Vietoris sequence 
for the spaces $\Sspace_t$, $\Uspace_t$, $\Vspace_t$ and $\Mspace_t$:
\begin{eqnarray*}
  \ldots  \to  \Hgroup_{p+1} (\Sspace_t)
          \to  \Hgroup_p (\Mspace_t)              
          \to  \Hgroup_p (\Uspace_t) \oplus \Hgroup_p (\Vspace_t) \\
          \to  \Hgroup_p (\Sspace_t)
          \to  \ldots .
\end{eqnarray*}
For $0 \leq t < 1$, $\Sspace_t$ is a point or
a closed $(n+1)$-dimensional ball,
and its only non-trivial homology group is
of dimension $0$, with rank $\Betti_0 (\Sspace_t) = 1$.
It follows that for $1 \leq p \leq n$,
the groups defined by $\Mspace_t$, $\Uspace_t$, $\Vspace_t$
are surrounded by trivial groups,
which implies that the groups of $\Mspace_t$
and the direct sums of the groups
of $\Uspace_t$ and $\Vspace_t$ are isomorphic.
We thus get
\begin{eqnarray}
  \Betti_0 (\Mspace_t) &=& \Betti_0 (\Uspace_t) + \Betti_0 (\Vspace_t) - 1,
  \label{eqn:MVBall-0}                                                       \\
  \Betti_p (\Mspace_t) &=& \Betti_p (\Uspace_t) + \Betti_p (\Vspace_t),
  \label{eqn:MVBall}
\end{eqnarray}
for $1 \leq p \leq n$.
We may also consider the relative homology groups,
again connected by a Mayer-Vietoris sequence:
\begin{eqnarray*}
  \ldots  \to  \Hgroup_{p+1} (\Sspace, \Sspace^t)
          \to  \Hgroup_p (\Mspace, \Mspace^t)           \\
          \to  \Hgroup_p (\Uspace, \Uspace^t) \oplus
               \Hgroup_p (\Vspace, \Vspace^t)           
          \to  \Hgroup_p (\Sspace, \Sspace^t)
          \to  \ldots .
\end{eqnarray*}
For $0 < t \leq 1$, the only non-trivial homology group
of the pair $(\Sspace, \Sspace^t)$ is in dimension $n+1$, which implies
\begin{eqnarray}
  \Betti_p (\Mspace, \Mspace^t)  &=&  \Betti_p (\Uspace, \Uspace^t)
                                    + \Betti_p (\Vspace, \Vspace^t),
  \label{eqn:MVopenBall} \\
  \Betti_n (\Mspace, \Mspace^t)  &=&  \Betti_n (\Uspace, \Uspace^t)
                                    + \Betti_n (\Vspace, \Vspace^t) + 1,
  \label{eqn:MVopenBall-n}
\end{eqnarray}
for $0 \leq p \leq n-1$.

\paragraph{Separating manifold.}
Combining Alexander duality with Mayer-Vietoris, we get relations
between the homology of the sublevel sets of $\Mspace$ and $\Uspace$.
More specifically, we get
\begin{eqnarray}
  \Betti_p (\Mspace)             &=&  \Betti_p (\Uspace)
                                    + \Betti_q (\Uspace) ,      
  \label{eqn:MVADSphere}                                            \\
  \Betti_p (\Mspace_t)           &=&  \Betti_p (\Uspace_t)
                                    + \Betti_q (\Uspace, \Uspace^t) ,      
  \label{eqn:MVADBall}                                              \\
  \Betti_p (\Mspace, \Mspace^t)  &=&  \Betti_p (\Uspace, \Uspace^t)
                                    + \Betti_q (\Uspace_t) ,
  \label{eqn:MVADopenBall}
\end{eqnarray}
for all $0 \leq p \leq n$.
Here, we combine
\eqref{eqn:MVSphere-0}, \eqref{eqn:MVSphere}, \eqref{eqn:MVSphere-n} with
\eqref{eqn:ADSphere-0}, \eqref{eqn:ADSphere}, \eqref{eqn:ADSphere-n}
to get \eqref{eqn:MVADSphere}.
Similarly, we combine \eqref{eqn:MVBall-0}, \eqref{eqn:MVBall} with
\eqref{eqn:ADBall-0}, \eqref{eqn:ADBall} to get \eqref{eqn:MVADBall}.
Finally, we exploit the symmetry between $\Uspace$ and $\Vspace$ and combine
\eqref{eqn:MVopenBall}, \eqref{eqn:MVopenBall-n} with
\eqref{eqn:ADBall-0}, \eqref{eqn:ADBall}
to get \eqref{eqn:MVADopenBall}.
\ignore{
It follows that \eqref{eqn:MVADBall} and \eqref{eqn:MVADopenBall}
hold for $0 < t < 1$.
}
Note that \eqref{eqn:MVADBall} and \eqref{eqn:MVADopenBall}
imply $\Betti_p (\Mspace_t) = \Betti_q (\Mspace, \Mspace^t)$,
which is a consequence of Lefschetz duality; see \cite{Mun84}.

\begin{figure}[hbt]
 \vspace{0.0in}
 \centering
 \resizebox{!}{1.8in}{\input{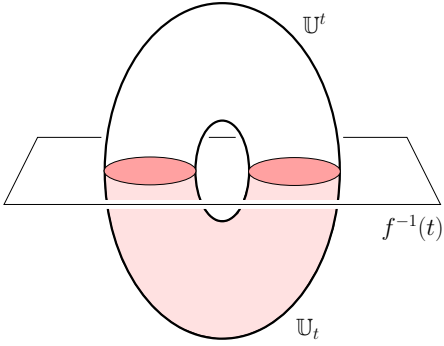}}
 \caption{The level set defined by $t$ splits the solid torus
   into two halves.}
 \label{fig:torus}
\end{figure}

\paragraph{An example.}
We illustrate the above relationships with an example.
Let $\Mspace$ be a $2$-dimensional torus in $\Sspace^3$.
Accordingly, $\Uspace$ and $\Vspace$ are the two solid tori
that decompose the $3$-sphere and intersect in $\Mspace$.
This is sketched in Figure \ref{fig:torus}, where we assume
that $\Uspace$ is the part of space surrounded by $\Mspace$,
while $\Vspace$ is the space surrounding $\Mspace$.
The only non-zero Betti numbers of the solid torus
are $\Betti_0 (\Uspace) = \Betti_1 (\Uspace) = 1$.
We thus get
\begin{eqnarray*}
  \Betti_0 (\Mspace)  &=&  \Betti_0 (\Uspace) + \Betti_2 (\Uspace)
              ~~~~~~~ ~~=~~  1 ,                                      \\
  \Betti_1 (\Mspace)  &=&  \Betti_1 (\Uspace) + \Betti_1 (\Uspace)
              ~~~~~~~ ~~=~~  2 ,                                      \\
  \Betti_2 (\Mspace)  &=&  \Betti_2 (\Uspace) + \Betti_0 (\Uspace) + 1
                      ~~=~~  1 
\end{eqnarray*}
from \eqref{eqn:MVADSphere}.
These are the correct Betti numbers of the $2$-dimensional torus.
Next, choose $t$ so that the sublevel set of $\Uspace$ is half
the solid torus.
Its only non-zero Betti number is $\Betti_0 (\Uspace_t) = 1$,
and the only non-zero relative Betti number of the pair is
$\Betti_1 (\Uspace, \Uspace^t) = 1$.
We thus get
\begin{eqnarray*}
  \Betti_0 (\Mspace_t)  &=&  \Betti_0 (\Uspace_t) + \Betti_2 (\Uspace, \Uspace^t)
                        ~~=~~  1 ,                                      \\
  \Betti_1 (\Mspace_t)  &=&  \Betti_1 (\Uspace_t) + \Betti_1 (\Uspace, \Uspace^t)
                        ~~=~~  1 ,                                      \\
  \Betti_2 (\Mspace_t)  &=&  \Betti_2 (\Uspace_t) + \Betti_0 (\Uspace, \Uspace^t)
                        ~~=~~  0 
\end{eqnarray*}
from \eqref{eqn:MVADBall}.
These are the correct Betti numbers of the sublevel set of $\Mspace$.
Finally, we get
\begin{eqnarray*}
  \Betti_0 (\Mspace, \Mspace^t)  &=&  \Betti_2 (\Uspace_t) + \Betti_0 (\Uspace, \Uspace^t)
                        ~~=~~  0 ,                                      \\
  \Betti_1 (\Mspace, \Mspace^t)  &=&  \Betti_1 (\Uspace_t) + \Betti_1 (\Uspace, \Uspace^t)
                        ~~=~~  1 ,                                      \\
  \Betti_2 (\Mspace, \Mspace^t)  &=&  \Betti_0 (\Uspace_t) + \Betti_2 (\Uspace, \Uspace^t)
                        ~~=~~  1 
\end{eqnarray*}
from \eqref{eqn:MVADopenBall}.
These are the correct relative Betti numbers of the pair $(\Mspace, \Mspace^t)$.

\section{Persistence}
\label{sec3}

Starting with a brief introduction of persistent homology,
we review its recent combinatorial expression as a point calculus;
see \cite{EdHa10} and \cite{BCE11,CSM09} for more comprehensive
discussions.

\paragraph{Background.}
We take the step from homology to persistent homology by replacing
a space with the sequence of (closed and open) sublevel sets
of a function on the space.
To explain this, let $\Xspace$ be compact
and $g: \Xspace \to [0,1]$ tame.
As before, we write $\Xspace_t = g^{-1} [0, t]$ for
the sublevel set and $\Xspace^t = g^{-1} [t, 1]$
for the superlevel set defined by $t$.
Since $g$ is tame, it has finitely many homological critical values,
$s_0$ to $s_m$.
We assume w.l.o.g.\ that $s_0=0$ and $s_m=1$,
else, we add them as additional critical values.
We interleave the $s_i$ with homological regular values $t_i$,
such that $s_0 < t_0 < s_1 < \ldots < t_{m-1} < s_{m}$.
Taking the homology of closed and open sublevel sets
at the regular values, we get the filtration
$$
  \Xgroup_0 \to \ldots \to \Xgroup_{m-1} \to \Xgroup_m \to \Xgroup_{m+1}
                       \to \ldots \to \Xgroup_{2m} ,
$$
where $\Xgroup_i$ is $\Hgroup (\Xspace_{t_i})$, if $0 \leq i \leq m-1$,
$\Hgroup (\Xspace)$, if $i = m$,
and $\Hgroup (\Xspace, \Xspace^{t_{2m-i}})$, if $m < i \leq 2m$.
For notational convenience, 
we add the trivial groups, $\Xgroup_{-1} = 0$ and
$\Xgroup_{2m+1} = 0$, at the beginning and end of the filtration.
The maps connecting the homology groups are induced by the
inclusions $\Xspace_s \subseteq \Xspace_t$ and
$\Xspace^t \subseteq \Xspace^s$, whenever $s \leq t$.
The maps compose and we write
$\gmap_{i,j} : \Xgroup_i \to \Xgroup_j$.
Reading the filtration from left to right, we see homology classes
appear and disappear.
To understand these events, we say a class $\alpha \in \Xgroup_i$
is \emph{born} at $\Xgroup_i$ if it does not belong to the image
of $\gmap_{i-1,i}$.
If furthermore $\gmap_{i,j} (\alpha)$ belongs to the image of $\gmap_{i-1,j}$
but $\gmap_{i,j-1} (\alpha)$ does not belong to
the image of $\gmap_{i-1,j-1}$,
then we say $\alpha$ \emph{dies entering} $\Xgroup_j$.
Since the filtration starts and ends with trivial groups, 
every homology class has well-defined birth and death values.
Every event is associated with the immediately preceding homological
critical value, namely with $s_{i}$ if the event is at $\Xgroup_i$,
and with $s_{m-i}$ if the event is at $\Xgroup_{m+i+1}$,
for $0 \leq i \leq m-1$.
Events at $\Xgroup_m$ are associated with $s_m=1$.
We represent a class by a dot in the plane
whose two coordinates mark its birth and its death.
More precisely, the coordinates signal the increase and decrease
of Betti numbers, and the dot represents an entire coset of classes
that are born and die with $\alpha$.
All dots have coordinates in $[0,1]$, by construction.

\begin{figure}[hbt]
 \vspace{0.0in}
 \centering
 \resizebox{!}{1.55in}{\input{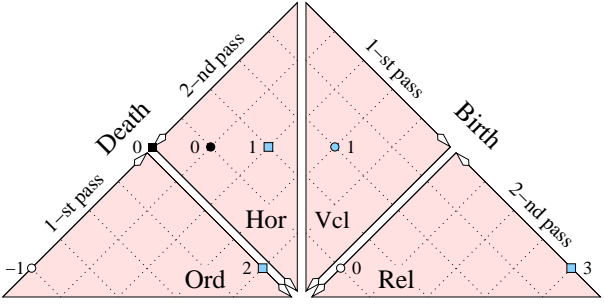}}
 \caption{The persistence diagram of the function restricted to the
   solid torus (dots drawn as circles),
   and to the complementary solid torus (dots drawn as squares),
   as sketched in Figure \protect{\ref{fig:torus}}.
   The white dots belong to the reduced diagram,
   the black dots belong to the standard diagram,
   and the shaded dots belong to both.}
 \label{fig:rdiagram}
\end{figure}

Note that we use each $s_i$ twice,
once during the first pass in which
the sublevel set grows from $\emptyset$ to $\Xspace$,
and then again during the second pass in which $\Xspace$ minus the
superlevel set shrinks back from $\Xspace$ to $\emptyset$.
When we collect the dots, we separate the passes
by drawing each coordinate axis twice.
The result is a multiset of dots, which we refer to as the
\emph{persistence diagram} of the function,
and denote as $\Ddgm (g)$, or $\Ddgm_p (g)$ if we restrict ourselves
to the dots representing $p$-dimensional classes;
see Figure \ref{fig:rdiagram}, where we label each dot with the
dimension of the classes it represents.
We further distinguish four regions within the diagram:
the \emph{ordinary}, the \emph{horizontal}, the \emph{vertical},
and the \emph{relative subdiagrams},
denoted as $\Odgm (g)$, $\Hdgm (g)$, $\Vdgm (g)$, and $\Rdgm (g)$.
For example, a dot belongs to the ordinary subdiagram
if its birth and death both happen during the first pass;
see again Figure \ref{fig:rdiagram}.

\paragraph{Reduced persistence diagrams.}
It is sometimes convenient to use reduced instead of standard
homology groups in the filtration.
There are small differences caused by the $(-1)$-dimensional class,
which exists if the space is empty.
As a first step, we introduce reduced versions of the relative homology groups,
which are isomorphic to the standard relative groups, for all $p$.
To define them, let $\omega$ be a new (dummy) vertex,
write $\omega \cdot \Xspace^t$ for the cone of $\omega$ over $\Xspace^t$,
and let $\Xspace_\omega^t = \Xspace \cupsp \omega \cdot \Xspace^t$
be the result of gluing the cone along $\Xspace^t$ to $\Xspace$.
Then the \emph{reduced relative homology group} of the pair
is $\rHgroup_p (\Xspace, \Xspace^t) = \rHgroup_p (\Xspace_\omega^t)$.
For example, if $\Xspace^t = \emptyset$,
then $\omega$ forms a separate component,
so that the reduced Betti number is
equal to the number of components of $\Xspace$,
just as $\Betti_0 (\Xspace, \Xspace^t)$.
The \emph{reduced persistence diagram} is now defined as before,
but for the filtration of reduced homology groups.
To describe this, we suppress the dimension and write
$\rXgroup_i$ for $\rHgroup (\Xspace_{t_i})$, if $0\leq i < m$, 
for $\rHgroup(\Xspace)$, if $i = m$, and
for $\rHgroup (\Xspace_\omega^{t_{2m-i}})$, if $m < i \leq 2m$.
The resulting sequence of reduced homology groups,
from $\rXgroup_0$ to $\rXgroup_{2m}$,
is connected from left to right by maps induced by inclusion.
Finally, we define $\rDdgm (g)$ by matching up the births and the deaths
and by drawing each pair as a dot in the plane.
Similar to before, we write $\rDdgm_p (g)$ when we restrict
ourselves to homology classes of dimension $p$.

For dimension $p \geq 1$, the reduced diagrams are the same
as the standard ones, simply because the groups are the same.
More formally, the persistence diagrams are the same
because the two filtrations form a commutative diagram
with vertical isomorphisms:
$$
  \begin{array}{ccccccccccc}
    0 & \to & \Xgroup_0^p  & \to & \Xgroup_1^p  &  \to & \ldots & \to & \Xgroup_{2m}^p &  \to & 0 \\

  \downarrow & &   \downarrow   &     & \downarrow   &      &        &     & \downarrow & & \downarrow         \\
  0 & \to & \rXgroup_0^p & \to & \rXgroup_1^p &  \to & \ldots & \to & \rXgroup_{2m}^p & \to & 0 ;
  \end{array}
$$
see the Persistence Equivalence Theorem \cite[page 159]{EdHa10}.
This is not true for $p = -1$, where the standard diagram is empty,
while the reduced diagram contains a single dot marking the transition
from an empty to a non-empty space.
To describe the difference for $p = 0$, we call a dot $(u,w)$
in the $0$-th standard persistence diagram \emph{extreme}
if no other dot has a first coordinate smaller than $u$ and a second
coordinate larger than $w$.
Here, we assume for simplicity that no two dots in the $0$-th diagram
share the same first coordinate or the same second coordinate.

We note that only dots in the horizontal subdiagram can be extreme.
The non-extreme dots also belong to the
reduced diagram of $g$, while the extreme dots exchange their coordinates
to form new dots in the reduced diagram; see Figure \ref{fig:rdiagram}.
The way the coordinates are exchanged will be important later,
so we describe the details of this process now.

\begin{figure}[hbt]
 \vspace{0.0in}
 \centering
 \resizebox{!}{1.55in}{\input{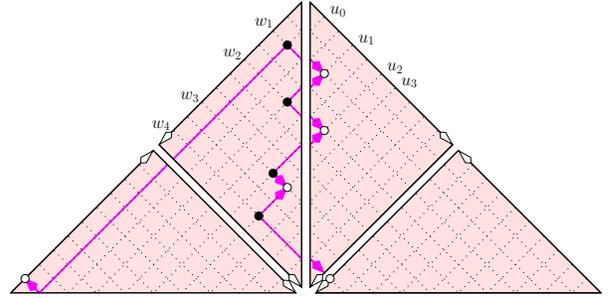}}
 \caption{The cascade combines the coordinates of the four black circles
   with $0$ and $1$ to form the five white circles.}
 \label{fig:cascade}
\end{figure}
\paragraph{Cascades.}
We let $\ell+1$ be the number of extreme dots and,
for a reason that will become clear shortly,
denote their coordinates with different indices
as $(u_k, w_{k+1})$, for $0 \leq k \leq \ell$,
ordering them such that $u_0 < u_1 < \ldots < u_\ell$.
By definition of extreme,
the corresponding second coordinates satisfy
$w_1 < w_2 < \ldots < w_{\ell+1}$.
Observe that $u_0$ is the globally minimum value
and $w_{\ell+1}$ is the globally maximum value.
Let us now construct the reduced diagram, focusing on dimension $0$.
We start the upward pass with the birth of the $(-1)$-dimensional homology
class at $0$, which dies at $u_0$,
so we have $(0, u_0)$ in $\rDdgm_{-1} (g)$.
The minima of all other components mark the births of
$0$-dimensional classes.
They can be interpreted as gaps between components.
We start the downward pass with the birth of a $0$-dimensional class
at $1$, which dies at $w_{\ell+1}$,
so we have $(1, w_{\ell+1})$ in $\rDdgm_0 (g)$.
The maxima of all other components mark the deaths of $0$-dimensional classes.
If $w$ is the second coordinate of a non-extreme component,
then it marks the death of the gap that opened up when we passed
the corresponding first coordinate, $u$, during the upward pass.
Hence, $(u, w)$ belongs to $\rDdgm_0 (g)$, and we note that this dot
also belongs to the standard diagram.
However, if $w = w_k$ is the second coordinate of a extreme dot,
then it marks the death of the gap that has opened up when we
passed the minimum of the next component with extreme dot
during the upward pass.
This minimum value is $u_k$,
so $(u_k, w_k)$ belongs to $\rDdgm_0 (g)$, for $1 \leq k \leq \ell$.

In summary, we see that contiguous extreme dots exchange one
pair of coordinates, using $0$ on the left and $1$ on the right
to complete the process.
More precisely, the reduced diagram is the same as the standard one,
except for substituting the dots \eqref{eqn:cascade-new}
for the dots \eqref{eqn:cascade-old}:
\begin{eqnarray}
  (u_0, w_1), (u_1, w_2), \ldots, (u_{\ell-1}, w_\ell), (u_\ell, w_{\ell+1}) ,
    \label{eqn:cascade-old} \\
  (0, u_0), (u_1, w_1), (u_2, w_2), \ldots, (u_\ell, w_\ell), (1, w_{\ell+1}) .
    \label{eqn:cascade-new}
\end{eqnarray}
For later use, we call this substitution a \emph{cascade}
and write $\rDdgm (g) = \Ddgm (g)^C$; see Figure \ref{fig:cascade}.
Note that all dots belong to the horizontal and vertical subdiagrams of dimension $0$,
except for $(1, w_{\ell+1})$, which belongs to the relative subdiagram
of dimension $0$,
and $(0, u_0)$, which belongs to the ordinary subdiagram of dimension $-1$.
We will often have $u_0 = 0$ and $w_{\ell+1} = 1$, in which case the
first and the last dots in \eqref{eqn:cascade-new} lie on the diagonal
and can be ignored.

\paragraph{Point calculus.}
As described in \cite{BCE11,CSM09}, the persistence diagram can be harvested
for a wealth of homological information.
This includes the ranks of the homology groups of the sets
$\Xspace_t$ and of the pairs $(\Xspace, \Xspace^t)$.
\begin{figure}[hbt]
 \vspace{0.0in}
 \centering
 \resizebox{!}{1.7in}{\input{Figs/diagram.pstex_t}}
 \caption{The persistence diagram of the function
   $g = f|_\Mspace$, as defined in Figure \protect{\ref{fig:torus}};
   the circle dots also belong to the diagram of $f|_\Uspace$,
   but the square dots do not.
   The rectangle on the left, $\Lwing_t (g)$, represents the homology
   of $\Mspace_t$, while the rectangle on the right, $\Rwing_t (g)$,
   represents the homology of $(\Mspace, \Mspace^t)$;
   compare with the Betti numbers computed for the example
   at the end of Section \protect{\ref{sec2}.}}
 \label{fig:diagram}
\end{figure}
To explain this, we write $\Lwing_t^p (g)$ for the multiset
of dots in the rectangle with lower corner $(t,t)$ on the base
edge of $\Odgm_p (g)$.
Similarly, we write $\Rwing_t^p (g)$ for the multiset of dots
if the lower corner, $(t,t)$, lies on the base edge of $\Rdgm_p (g)$;
see Figure \ref{fig:diagram}.
Often, we drop the dimension from the notation and write
$\Lwing_t (g)$ for the disjoint union of the multisets
$\Lwing_t^p (g)$, over all $p$, and similar for $\Rwing_t (g)$.
We can read the Betti numbers of $\Xspace_t$ by collecting
the dots in $\Lwing_t (g)$,
and we can read the relative Betti numbers of
$(\Xspace, \Xspace^t)$ by collecting the dots in $\Rwing_t (g)$.
Writing $\rLwing_t^p (g)$ and $\rRwing_t^p (g)$ for the corresponding
multisets in the reduced diagrams,
we note that the same relations hold between the reduced Betti numbers
and the reduced persistence diagrams.

\paragraph{Betti numbers and persistence diagrams.}
To motivate the first result of this paper, we rewrite
\eqref{eqn:ADSphere-0} to \eqref{eqn:ADBall} for reduced Betti numbers:
\begin{eqnarray}
  \rBetti_p (\Vspace_t)  &=&  \Betti_q (\Uspace, \Uspace^t) ,
  \label{eqn:reducedAD}
\end{eqnarray}
for all $0 \leq p \leq n$ and all regular values $t$.
Observe also that $\Vspace_t = \emptyset$ iff
$\Uspace - \Uspace^t = \Sspace - \Sspace^t$, which implies that
\eqref{eqn:reducedAD} also holds for $p = -1$.
In words, for every dimension $p$ and every regular value $t$, 
the rectangles $\rLwing_t^p (f|_\Vspace)$ and $\rRwing_t^q (f|_\Uspace)$
have equally many dots.
This $1$-parameter family of relations is satisfied if the
reduced diagrams of $f|_\Vspace$ and $f|_\Uspace$ are
reflections of each other.
This is our first result, stated as the Land and Water Theorem in
Section \ref{sec4}.
We see an illustration in Figure \ref{fig:rdiagram},
which shows the reduced and non-reduced diagrams of the function
$f$ in Figure \ref{fig:torus} restricted to the solid torus, $\Uspace$,
and to the complementary solid torus, $\Vspace$.
Removing the black dots, we are left with the two reduced diagrams,
which are indeed reflections of each other.

To motivate our second result, we consider the relations
\eqref{eqn:MVADSphere} to \eqref{eqn:MVADopenBall}.
They say that for every regular value $t$,
the rectangle $\Lwing_t^p (f|_\Mspace)$
has the same number of dots as
$\Rwing_t^q (f|_\Uspace)$ and $\Lwing_t^p (f|_\Uspace)$ together.
Similarly, $\Rwing_t^p (f|_\Mspace)$
has the same number of dots as $\Lwing_t^q (f|_\Uspace)$
and $\Rwing_t^p (f|_\Uspace)$ together.
This $1$-parameter family of relations is satisfied if the
persistence diagram of $f|_\Mspace$ is the disjoint union
of the diagram of $f|_\Uspace$ and of its reflection.
This is our second result, stated as the Euclidean Shore Theorem
in Section \ref{sec5}.
It is illustrated in Figure \ref{fig:diagram}, which shows
the persistence diagram of the function
$f|_\Mspace$ in Figure \ref{fig:torus}.
Comparing Figures \ref{fig:rdiagram} and \ref{fig:diagram},
we see that the two circle dots in Figure \ref{fig:diagram}
also belong to the diagram of $f|_\Uspace$,
whereas the two square dots are reflections of the circles.

While \eqref{eqn:ADSphere-0} to \eqref{eqn:ADBall}
and \eqref{eqn:MVADSphere} to \eqref{eqn:MVADopenBall}
motivate our two results, these relations are not sufficient to prove them.
Indeed, \eqref{eqn:MVADSphere} to \eqref{eqn:MVADopenBall} hold in general,
but the Euclidean Shore Theorem requires that the minimum and maximum of the
perfect Morse function belong to a common component of $\Vspace$.
If this condition is violated, then we can have coordinate exchanges
among the dots that contradict the Euclidean Shore Theorem
without affecting the relations \eqref{eqn:MVADSphere} to \eqref{eqn:MVADopenBall}.
An example of this phenomenon can be seen in Figure \ref{fig:counterexample},
which will be discussed in Section \ref{sec5}.

\section{Land and Water}
\label{sec4}

In this section, we present our first result,
which extends Alexander duality from spaces to functions.
Its proof needs a general algebraic result about
contravariant filtrations, which we present first.

\paragraph{Compatible pairings.}
A \emph{pairing} between two finite-di\-mensional vector spaces
$\Xgroup$ and $\Ygroup$ over the field $\Fgroup$ is a bilinear map
$\scalprod{}{} : \Xgroup \times \Ygroup \to \Fgroup$.
The pairing is \emph{non-degenerate} if every $\xi \in \Xgroup$
has at least one $\eta \in \Ygroup$ with $\scalprod{\xi}{\eta} = 1$,
and, symmetrically, every element $\eta \in \Ygroup$ has at least one
$\xi \in \Xgroup$ with $\scalprod{\xi}{\eta} = 1$.
The term is justified by the fact that a non-degenerate pairing
implies an isomorphism between the vector spaces.
To see this, we choose arbitrary bases of $\Xgroup$ and $\Ygroup$,
writing $n_\Xgroup$ and $n_\Ygroup$ for their sizes,
and represent the pairing by its values between the basis elements.
This gives an invertible $0$-$1$ matrix with $n_\Xgroup$ rows and $n_\Ygroup$ columns.
Reducing the matrix to diagonal form gives the claimed isomorphism.

\begin{figure}[hbt]
 \vspace{0.0in}
 \centering
 \resizebox{!}{2.3in}{\input{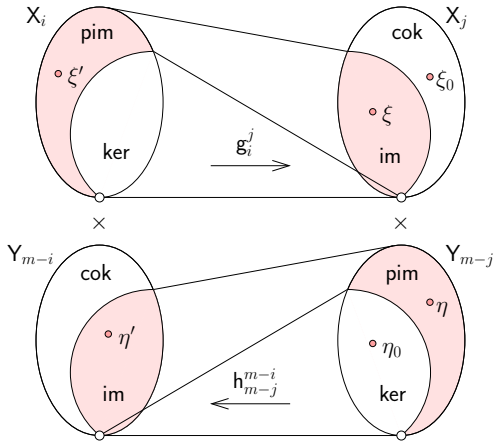}}
 \caption{Splitting the spaces in the two
   anti-parallel filtrations.}
 \label{fig:split}
\end{figure}

Now consider two filtrations, each consisting of $m+1$
finite-dimensional vector spaces over $\Fgroup$,
and non-degenerate pairings
connecting the filtrations contravariantly:
$$
  \begin{array}{ccccccc}
   \Xgroup_0  & \to   & \Xgroup_1     & \to   & \ldots & \to   & \Xgroup_m   \\
   \times     &       & \times        &       &        &       & \times      \\
   \Ygroup_m  & \from & \Ygroup_{m-1} & \from & \ldots & \from & \Ygroup_0 .
  \end{array}
$$
It will be convenient to assume that the filtrations begin and end
with the trivial vector space:
$\Xgroup_0 = \Xgroup_m = \Ygroup_0 = \Ygroup_m = 0$.
Write $\gmap_i^j : \Xgroup_i \to \Xgroup_j$ and
$\hmap_i^j : \Ygroup_i \to \Ygroup_j$ for the maps
upstairs and downstairs.
We call the pairings \emph{compatible} with these maps
if $\scalprod{\xi}{\eta} = \scalprod{\xi'}{\eta'}$
for every $\xi' \in \Xgroup_i$ and $\eta \in \Ygroup_{m-j}$,
where $\xi = \gmap_i^j (\xi')$ and $\eta' = \hmap_{m-j}^{m-i} (\eta)$.
By what we said above, we have isomorphisms
$\isomap_j : \Xgroup_j \to \Ygroup_{m-j}$, for all $0 \leq j \leq m$,
but we need more, namely isomorphisms that satisfy a compatibility
condition with the horizontal maps.
Before defining what we mean, we prove a few properties of
the compatible pairings.
To begin, we split the vector spaces connected by $\gmap_i^j$:
\begin{eqnarray}
  \Xgroup_i  &=&  \pimage{\gmap_i^j} \oplus \kernel{\gmap_i^j} , \\
  \Xgroup_j  &=&  \image{\gmap_i^j} \oplus \cokernel{\gmap_i^j} ,
\end{eqnarray}
and similarly for $\Ygroup_{m-j}$ and $\Ygroup_{m-i}$;
see Figure \ref{fig:split}.
Here, we have of course
$\pimage{\gmap_i^j} = \Xgroup_i / \kernel{\gmap_i^j}$
and $\cokernel{\gmap_i^j} = \Xgroup_j / \image{\gmap_i^j}$.
The restriction of the pairing to the kernel of $\gmap_i^j$ times
the image of $\hmap_{m-j}^{m-i}$ is well-defined because elements
of the kernel and the image are also elements of the vector spaces.
In contrast, an element of a quotient is an entire coset,
whose members may have inconsistent values under the pairing.
To define the restriction, we assume a basis that splits along the same lines,
so that the coset in $\pimage{\gmap_i^j}$ has a unique representation
of the form $\xi' + \kernel{\gmap_i^j}$,
in which the projection of $\xi'$ to the kernel is zero.
Finally, we define the \emph{restricted pairing} using $\xi'$
as a proxy for the coset it defines.
Similarly, we define the restriction of the pairing to the
cokernel of $\gmap_i^j$.
We have the following properties:
\begin{enumerate}
  \item[(A)]  The pairing restricted to
    $\image{\gmap_i^j} \times \kernel{\hmap_{m-j}^{m-i}}$ is trivial,
    that is: all pairs map to $0$.
  \item[(B)]  We can find a basis of the cokernel such that
    the pairing restricted to
    $\cokernel{\gmap_i^j} \times \pimage{\hmap_{m-j}^{m-i}}$
    is trivial.
\end{enumerate}
To prove (A), let $\xi \in \image{\gmap_i^j}$ and
$\eta_0 \in \kernel{\hmap_{m-j}^{m-i}}$ both be non-zero.
The image of $\eta_0$ is $0$, and a preimage,
$\xi'$, of $\xi$ exists.
This implies $\scalprod{\xi}{\eta_0} = \scalprod{\xi'}{0} = 0$,
as claimed.

To prove (B),
let $\xi_0 \in X_j$ and $\eta \in Y_{m-j}$ both be non-zero,
with zero projection of $\xi_0$ to $\image{\gmap_i^j}$
and of $\eta$ to $\kernel{\hmap_{m-j}^{m-i}}$.
Recall that $\xi_0$ and $\eta$ act as proxies for the cosets
$\xi_0 + \image{\gmap_i^j}$ in $\cokernel{\gmap_i^j}$ and
$\eta + \kernel{\hmap_{m-j}^{m-i}}$ in $\pimage{\hmap_{m-j}^{m-i}}$.
If $\scalprod{\xi_0}{\eta} = 0$, then there is nothing to do,
else we take steps to set the value
of the pairing for the cosets to zero.
Assume without loss of generality that $\xi_0$ is
an element of the chosen basis of $X_j$.
By the non-degeneracy of the pairing $\Xgroup_i \times \Ygroup_{m-i}$,
there exists $\xi' \in X_j$ with $\scalprod{\xi'}{\eta'} = 1$,
where $\eta' = \hmap_{m-j}^{m-i} (\eta)$.
By (A), we can assume that the projection of $\xi'$ to
the $\kernel{\gmap_i^j}$ is zero.
By compatibility, we have $\scalprod{\xi}{\eta} = 1$,
where $\xi = \gmap_i^j (\xi')$.
Now we have $\scalprod{\xi_0}{\eta} = \scalprod{\xi}{\eta} = 1$,
and $\xi_0 \neq \xi$ because one belongs to the cokernel and
the other to the image of $\gmap_i^j$.
Replacing $\xi_0$ by $\xi_1 = \xi_0 - \xi$ in the basis of $X_j$,
we get $\scalprod{\xi_1}{\eta} = 0$ by bilinearity.
We note that we still have a basis that splits as required.
The operation achieves its goal by changing the restricted pairing
without affecting the pairing.
Repeating the operation for other basis elements of the cokernel,
we eventually get a basis such that the pairing restricted to
$\cokernel{\gmap_i^j} \times \pimage{\hmap_{m-j}^{m-i}}$ is trivial,
as claimed.

\begin{figure}[hbt]
 \vspace{0.0in}
 \centering
 \resizebox{!}{1.0in}{\input{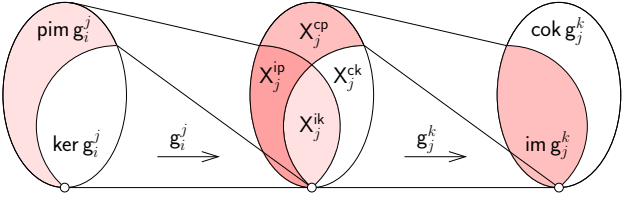}}
 \caption{Splitting the vector spaces into four.}
 \label{fig:four-way}
\end{figure}

\paragraph{Compatible isomorphisms.}
We are now ready to specify the compatibility condition for the
vertical isomorphisms.
First, we require that they split along the same lines as the spaces.
Specifically, we consider indices $0 \leq i < j < k \leq m$ and the
maps $\gmap_i^j$ and $\gmap_j^k$, and we define
\begin{eqnarray}
  \Xgroup_j^{\ik}  &=&  \image{\gmap_i^j} \capsp \kernel{\gmap_j^k},  \\
  \Xgroup_j^{\ip}  &=&  \image{\gmap_i^j} / \Xgroup_j^{\ik},           \\
  \Xgroup_j^{\ck}  &=&  \kernel{\gmap_j^k} / \Xgroup_j^{\ik},          \\
  \Xgroup_j^{\cp}  &=&  \Xgroup_j / (\Xgroup_j^{\ik} \oplus
                        \Xgroup_j^{\ip} \oplus \Xgroup_j^{\ck}) ;
\end{eqnarray}
see Figure \ref{fig:four-way}.
For the case $i+1 = j = k-1$, these spaces represent the classes
that end at $\Xgroup_j$, that extend before and after $\Xgroup_j$,
that exist only at $\Xgroup_j$, and that begin at $\Xgroup_j$,
in this sequence.
Similarly, we consider the anti-parallel maps,
$\hmap_{m-k}^{m-j}$ and $\hmap_{m-j}^{m-i}$, and
we define the subspaces $\Ygroup_{m-j}^{\ik}$, $\Ygroup_{m-j}^{\ip}$,
$\Ygroup_{m-j}^{\ck}$, $\Ygroup_{m-j}^{\cp}$ downstairs.
We call the vertical isomorphisms \emph{compatible} with the horizontal maps
if $\isomap_j: \Xgroup_j \to \Ygroup_{m-j}$ splits into isomorphisms
\begin{eqnarray}
  \isomap_j^{\ik}  &=&  \Xgroup_j^{\ik} \to \Ygroup_{m-j}^{\cp} ,           \\
  \isomap_j^{\ip}  &=&  \Xgroup_j^{\ip} \to \Ygroup_{m-j}^{\ip} ,           \\
  \isomap_j^{\ck}  &=&  \Xgroup_j^{\ck} \to \Ygroup_{m-j}^{\ck} ,           \\
  \isomap_j^{\cp}  &=&  \Xgroup_j^{\cp} \to \Ygroup_{m-j}^{\ik} ,
\end{eqnarray}
and the maps between the subspaces compose to the identity,
for every $0 \leq j \leq m$; see Figure \ref{fig:split}.
More formally, this means that the maps satisfy
\begin{eqnarray}
  \xi  &=&  \gmap_i^j (\isomap_i^{-1} (\hmap_{m-j}^{m-i} (\isomap_j (\xi)))) ,
  \label{eqn:identity}
\end{eqnarray}
for all $\xi \in \image{g_i^j}$ and $0 \leq i < j \leq m$.

We now prove that the compatible pairings imply the existence of
compatible isomorphisms.
The proof is constructive.
To start, we choose bases of $\Xgroup_j$ and $\Ygroup_{m-j}$ that split
four ways, as do the spaces.
Representing the pairing between $\Xgroup_j$ and $\Ygroup_{m-j}$ by the
thus defined $0$-$1$ matrix, we get a decomposition into $16$ blocks;
see Figure \ref{fig:blocks}.
\begin{figure}[hbt]
 \vspace{0.0in}
 \centering
 \resizebox{!}{1.4in}{\input{Figs/blocks.pstex_t}}
 \caption{The pairing between $\Xgroup_j$ and $\Ygroup_{m-j}$
   splits into four non-degenerate restrictions.}
 \label{fig:blocks}
\end{figure}
We use (A) and (B) to show that the bases can be chosen so that
all non-diagonal blocks are zero.
Indeed, the upper-right four blocks are zero because
$\Xgroup_j^{\ik} \oplus \Xgroup_j^{\ip} = \image{\gmap_i^j}$ and
$\Ygroup_{m-j}^{\ck} \oplus \Ygroup_{m-j}^{\ik}
  = \kernel{\hmap_{m-j}^{m-i}}$.
The lower-left four blocks are zero because
$\Xgroup_j^{\ck} \oplus \Xgroup_j^{\cp} = \cokernel{\gmap_i^j}$ and
$\Ygroup_{m-j}^{\cp} \oplus \Ygroup_{m-j}^{\ip}
  = \pimage{\hmap_{m-j}^{m-i}}$.
Repeating the argument with reversed roles
(swapping $X_j$ with $Y_{m-j}$), we conclude that the
remaining four non-diagonal blocks are also zero.
Since the entire matrix is invertible, this implies that all
four diagonal blocks are invertible.
We finally get the isomorphisms by turning the diagonal blocks into
unit matrices.
All we need here are standard row and column operations,
which preserve the pairing while changing
the basis for its description.
The compatibility condition \eqref{eqn:identity} is therefore but a
reformulation of the compatibility condition for the pairings.
This completes the argument for the existence of compatible isomorphisms.
The reason for our interest is the following straightforward result.
\begin{result}[Contravariant PE Theorem]
  Let $\Xgroup_0$ to $\Xgroup_m$ and $\Ygroup_0$ to $\Ygroup_m$
  be two filtrations contravariantly connected by compatible isomorphisms.
  Then there is a bijection between the classes upstairs and downstairs
  such that a class is born at $\Xgroup_i$ and dies entering $\Xgroup_j$
  iff the corresponding class is born at $\Ygroup_{m-j+1}$
  and dies entering $\Ygroup_{m-i+1}$.
\end{result}
\proof
 Recall that $\isomap_j: \Xgroup_j \to \Ygroup_{m-j}$ splits into four
 isomorphisms, for each choice of triple indices $i < j < k$.
 The case $i+1 = j = k-1$ implies that the births and deaths upstairs
 correspond to the deaths and births downstairs.
 Condition \eqref{eqn:identity} implies that the births and deaths
 are matched upstairs the same way as downstairs, only backward.
\eop

We see that the compatibility condition \eqref{eqn:identity}
is for contravariantly connected filtrations what commutativity
is for covariantly connected filtrations.
While the (covariant) Persistence Equivalence Theorem \cite[page 159]{EdHa10}
implies equal persistence diagrams,
the Contravariant P(ersistence) E(quivalence) Theorem
implies reflected persistence diagrams.

\paragraph{Alexander duality for functions.}
After getting the abstract prerequisites settled,
we are ready to state and prove our first result.
Let $f: \Sspace \to [0,1]$ be a perfect Morse function
on the $(n+1)$-dimensional sphere, and
$\Sspace = \Uspace \cupsp \Vspace$ a decomposition into
complementary subsets whose intersection, $\Mspace = \Uspace \capsp \Vspace$,
is an $n$-manifold.
Assuming $f|_\Mspace$ is tame, we have only finitely many
homological critical values, including $0$ and $1$.
Letting these critical values be $s_0$ to $s_m$, we interleave them
with homological regular values to get
$s_0 < t_0 < s_1 < \ldots < t_{m-1} < s_m$.
We are interested in the filtrations defined by the restrictions of
$f$ to $\Uspace$ and to $\Vspace$.
To describe them, we write $\rUgroup_i^p$ for the $p$-th reduced homology
group of $\Uspace_{t_i}$, if $0 \leq i < m$,
of $\Uspace$, if $i = m$,
and of $(\Uspace, \Uspace^{t_{2m-i}})$, if $m < i \leq 2m$.
In the same way, we define $\rVgroup_i^q$ for $0\leq i\leq 2m$
with respect to $\Vspace$.
Using Alexander duality to form isomorphisms, we get two contravariantly
connected filtrations:
$$
  \begin{array}{ccccccccc}
    \rUgroup_0^p      & \to   & \ldots & \to   & \rUgroup_m^p     & \to   & \ldots & \to   & \rUgroup_{2m}^p   \\
    \downarrow        &       &        &       & \downarrow       &       &        &       & \downarrow          \\
    \rVgroup_{2m}^q & \from & \ldots & \from & \rVgroup_{m}^q & \from & \ldots & \from & \rVgroup_0^q ,
  \end{array}
$$
where $p + q = n$.
We will prove shortly that the isomorphisms implied by Alexander
duality are compatible with the horizontal maps.
Assuming this much, we get our first result, which we formulate
using the superscript `$T$' for the operation that
reflects a dot across the vertical axis of a persistence diagram
and, at the same time, changes its dimension from $p$ to $q = n-p$.
\begin{result}[Land and Water Theorem]
  Let $\Uspace$ and $\Vspace$ be complementary subsets of
  $\Sspace = \Sspace^{n+1}$,
  and let $f: \Sspace \to [0,1]$ be a perfect Morse function
  whose restriction to the $n$-manifold $\Mspace = \Uspace \capsp \Vspace$ is tame.
  Then $\rDdgm (f|_\Vspace) = \rDdgm (f|_\Uspace)^T$.
\end{result}
\proof
 By the Contravariant PE Theorem, we have a bijection between the classes
 upstairs and downstairs that respects the matching of births and deaths.
 Specifically, if a class $\alpha$ is born at $\rUgroup_i^p$ and dies
 entering $\rUgroup_j^p$ then its corresponding class is born at
 $\rVgroup_{2m-j+1}^q$ and dies entering $\rVgroup_{2m-i+1}^q$.
 If $1 \leq i < j \leq m$, then the class $\alpha$ is represented
 by $(s_i, s_j)$ in the $p$-th ordinary subdiagram of $f|_\Uspace$.
 The corresponding class is represented by $(s_j, s_i)$
 in the $q$-th relative subdiagram of $f|_\Vspace$.
 The reflection maps the first and second
 coordinate-axes of the ordinary subdiagram to the second and first
 coordinate-axes of the relative subdiagram.
 It follows that it maps $(s_i, s_j)$ in the former to $(s_j, s_i)$ in the latter subdiagram.
 Other cases are similar, and we conclude that
 $\rDdgm_p (f|_\Uspace)$ and $\rDdgm_q (f|_\Vspace)$
 are reflections of each other.
 Writing this more succinctly gives the claimed relationship
 between the reduced persistence diagrams of $f$ restricted
 to the two complementary subsets of the sphere.
\eop

\paragraph{Alexander pairing.}
We fill the gap in the proof of the Land and Water Theorem
by establishing compatible pairings between the groups
upstairs and downstairs.
For the sake of simplicity, we restrict ourselves
to the case $\Fgroup = \Zspace_2$.
While the argument for general fields is similar,
it requires oriented simplices, which is a technical formalism
we prefer to avoid.
We begin with the pairing implicit in
Lefschetz duality for manifolds with possibly non-empty boundary:
$$
  \scalprod{}{}_L: \rHgroup_p (\Uspace_t) \times \rHgroup_{q+1} (\Uspace_t, \partial\Uspace_t)
                   \to \Zspace_2 ,
$$
which is defined by mapping two classes to the parity of the
number of intersections between representing cycles.
We get such a pairing for every regular value $t$,
and these pairings are
compatible with the horizontal maps induced by inclusion of
the sublevel sets; see \cite{CEH09}.
Next, we reduce the dimension of the second factor using a mapping,
$\varphi$, which we compose from four simpler mappings:
\begin{eqnarray*}
  \rHgroup_{q+1} (\Uspace_t, \partial\Uspace_t)  &\stackrel{\varphi_1}{\to}&
           \rHgroup_{q+1} (\Sspace, \Sspace^t \cupsp \Vspace_t)       \\
                                       &\stackrel{\varphi_2}{\to}&
           \rHgroup_q (\Sspace^t \cupsp \Vspace_t)                    \\
                                       &\stackrel{\varphi_3}{\to}&
           \rHgroup_q (\Sspace^t \cupsp \Vspace_t, \Sspace^t)         \\
                                       &\stackrel{\varphi_4}{\to}&
           \rHgroup_q (\Vspace, \Vspace^t). 
\end{eqnarray*}
The first mapping, $\varphi_1$, is an isomorphism defined by excision.
Next, $\varphi_2$, is the connecting homomorphism of the
exact sequence of the pair $(\Sspace, \Sspace^t \cupsp \Vspace_t)$.
Since the $q$-th reduced homology of $\Sspace$ is trivial,
for all $q \neq n+1$, $\varphi_2$ is an isomorphism for all $0 \leq q < n$.
It is a surjection for $q = n$.
The third mapping, $\varphi_3$, is induced by inclusion.
It occurs in the exact sequence of
the pair $(\Sspace^t \cupsp \Vspace_t, \Sspace^t)$.
For $0 < t \leq 1$, the reduced homology groups of $\Sspace^t$ are all trivial,
implying that $\varphi_3$ is an isomorphism.
For $t = 0$, $\varphi_3$ is the trivial isomorphism.
Finally, $\varphi_4$ is again an isomorphism defined by excision.
In summary, we get
$\varphi: \rHgroup_{q+1} (\Uspace_t, \partial\Uspace_t) \to \rHgroup_q (\Vspace, \Vspace^t)$,
which is an isomorphism for $0 \leq q < n$ and a surjection for $q = n$.
More specifically, each class $\beta \in \rHgroup_n (\Vspace, \Vspace^t)$
has two preimages.
Indeed, an $n$-cycle representing the corresponding class in
$\rHgroup_n (\Sspace^t \cupsp \Vspace_t)$
partitions the components of $\Uspace_t$ into two subsets,
and each subset generates an $(n+1)$-dimensional relative homology
class that maps to $\beta$.

Note that the roles of $\Uspace$ and $\Vspace$ can be interchanged to get
a mapping from $\rHgroup_{q+1} (\Vspace_t, \partial\Vspace_t)$ to
$\rHgroup_q (\Uspace, \Uspace^t)$.
With this, we are ready to construct the \emph{Alexander pairing}:
$$
  \scalprod{}{}_A : \rHgroup_p (\Uspace_t) \times \rHgroup_q (\Vspace, \Vspace^t)
                      \to \Zspace_2 , 
$$
defined by $\scalprod{\alpha}{\beta}_A = \scalprod{\alpha}{\beta'}_L$,
where $\beta'$ is a preimage of $\beta$ under $\varphi$.
This fixes the pairing of $\rUgroup_i^p$ and $\rVgroup_{2m-i}^q$
for $0 \leq i < m$.
Similarly, we define $\scalprod{}{}_A$ for
$\rHgroup_p (\Uspace, \Uspace^t) \times \rHgroup_q (\Vspace_t)$ 
which fixes the pairing for $m <  i \leq 2m$.
The pairing is clearly well-defined for $0 < p < n$, where $\varphi$
is an isomorphism.
The remaining two cases are symmetric, and we consider $p = 0$
and $\beta \in \rHgroup_n (\Vspace, \Vspace^t)$.
As noted before, $\beta$ has two preimages, $\beta'$ and $\beta''$,
which are generated by complementary subsets of $\Uspace_t$.
Every $\alpha \in \rHgroup_0 (\Uspace_t)$ is represented by an even number
of points, which the complementary subsets partition into two even
or two odd subsets.
In either case, we have
$\scalprod{\alpha}{\beta'}_L = \scalprod{\alpha}{\beta''}_L$,
so the pairing is well-defined in all cases. 
Moreover, the compatibility
of the Alexander pairings within the left and right halves follows
from the compatibility of the Lefschetz pairings.

It remains to define the paring in the middle of the sequence,
for $\rUgroup_m^p = \rHgroup_p (\Uspace)$ and
$\rVgroup_m^q = \rHgroup_q (\Vspace)$.
To this end, let $\alpha$ be a reduced $p$-cycle in $\Uspace$,
$\beta$ a reduced $q$-cycle in $\Vspace$,
and $\beta'$ a $(q+1)$-chain whose boundary is $\beta$.
We define $\scalprod{\alpha}{\beta}_A$ 
by counting the intersections between $\alpha$ and $\beta'$.
Note that $\beta'$ is the preimage of $\beta$ under the previously
defined mapping,
$\varphi: \rHgroup_{q+1} (\Uspace_t, \partial\Uspace_t) \to \rHgroup_q (\Vspace, \Vspace^t)$, for any regular value $t$.
Therefore, the pairing in the middle is compatible with
the pairings in left half,
again exploiting the compatibility of the Lefschetz pairing.
Alternatively, we can define the pairing by taking
a $(p+1)$-chain $\alpha'$ with boundary $\alpha$ and counting
its intersections with $\beta$.
In this case, $\alpha'$ is the preimage of $\alpha$ under the map
$\varphi: \rHgroup_{p+1} (\Vspace_t, \partial\Vspace_t)
      \to \rHgroup_p (\Uspace, \Uspace^t)$,
and compatibility with the pairings in the right half follows.
Indeed, both definitions are equivalent,
as already observed by Lefschetz \cite{Lef26}.
We give a simple proof for the case $\Fgroup = \Zspace_2$:
\begin{result}[Bridge Lemma]
  Let $\alpha$ and $\beta$ be non-intersecting reduced cycles on the
  $(n+1)$-sphere whose dimensions add up to $n$.
  Let $\alpha'$ and $\beta'$ be chains whose boundaries are
  $\alpha$ and $\beta$, respectively.
  Then $\scalprod{\alpha'}{\beta}_L = \scalprod{\alpha}{\beta'}_L$.
\end{result}
\proof
 We can assume that $\alpha$ and $\beta$ do not intersect.
 The intersection of $\alpha'$ and $\beta'$ is a $1$-chain.
 Its boundary consists of an even number of points,
 and is the disjoint union of intersections of
 $\alpha'$ with $\beta$ and of $\alpha$ with $\beta'$.
 This implies that both types of intersections occur with the same parity.
\eop

\section{Shore}
\label{sec5}

This section presents our second result,
which extends the combination of Alexander duality
and Mayer-Vietoris sequences from spaces to functions.

\paragraph{Mayer-Vietoris sequence of filtrations.}
Assuming that $f|_\Mspace$ is tame,
we write $s_0 < t_0 < s_1 < \ldots < t_{m-1} < s_m$ for the
interleaved sequence of homological critical and regular values,
as before.
The main tool in this section is the diagram obtained
by connecting the filtrations of $f$ and its restrictions
with Mayer-Vietoris sequences.
We describe this using shorthand notation for the groups.
Consistent with earlier notation, we write $\Sgroup_i^p$
for the $p$-th homology group of $\Sspace_{t_i}$, for $0 \leq i < m$,
of $\Sspace$, for $i = m$,
and of $(\Sspace, \Sspace^{t_{2m-i}})$, for $m < i \leq 2m$.
Similarly, we write $\Mgroup_i^p$ for the $p$-th homology group
of $\Mspace_{t_i}$, of $\Mspace$, and of $(\Mspace, \Mspace^{t_{2m-i}})$.
Finally, we write $\Dgroup_i^p$ for the direct sum of the $p$-th
homology groups of $\Uspace_{t_i}$ and $\Vspace_{t_i}$,
of $\Uspace$ and $\Vspace$,
and of $(\Uspace, \Uspace^{t_{2m-i}})$ and $(\Vspace, \Vspace^{t_{2m-i}})$.
Drawing the filtrations from left to right and the Mayer-Vietoris
sequences from top to bottom, we get
$$
  \begin{array}{ccccccccccc}
    \downarrow       &     &        &     & \downarrow      &     &        &     & \downarrow             \\
    \Sgroup_0^{p+1}  & \to & \ldots & \to & \Sgroup_m^{p+1} & \to & \ldots & \to & \Sgroup_{2m}^{p+1}    \\
    \downarrow       &     &        &     & \downarrow      &     &        &     & \downarrow             \\
    \Mgroup_0^p      & \to & \ldots & \to & \Mgroup_m^p     & \to & \ldots & \to & \Mgroup_{2m}^p        \\
    \downarrow       &     &        &     & \downarrow      &     &        &     & \downarrow             \\
    \Dgroup_0^p      & \to & \ldots & \to & \Dgroup_m^p     & \to & \ldots & \to & \Dgroup_{2m}^p        \\
    \downarrow       &     &        &     & \downarrow      &     &        &     & \downarrow             \\
    \Sgroup_0^p      & \to & \ldots & \to & \Sgroup_m^p     & \to & \ldots & \to & \Sgroup_{2m}^p        \\
    \downarrow       &     &        &     & \downarrow      &     &        &     & \downarrow             
 \end{array}
$$
All squares commute, which is particularly easy to see for the squares
that connect groups of the same dimension,
whose maps are all induced by inclusion.
For $1 \leq p \leq n$, all groups $\Sgroup_i^p$ are trivial.
By exactness of the Mayer-Vietoris sequences, this implies that
the maps $\mmap_i^p: \Mgroup_i^p \to \Dgroup_i^p$
are isomorphisms, for $1 \leq p \leq n-1$.
The persistence diagram of $f|_\Mspace$ is therefore the disjoint union
of the persistence diagrams of $f|_\Uspace$ and $f|_\Vspace$.
More generally, we claim:
\begin{result}[General Shore Theorem]
  Let $n$ be a positive integer, 
  let $\Uspace$ and $\Vspace$ be complementary subsets of
  $\Sspace = \Sspace^{n+1}$,
  and let $f: \Sspace \to [0,1]$ be a perfect Morse function
  whose restriction to the $n$-manifold $\Mspace = \Uspace \capsp \Vspace$ is tame.
  Then
\begin{eqnarray}
  \Ddgm_0 (f|_\Mspace)  &=&  [ \Ddgm_0 (f|_\Uspace) \sqcup \Ddgm_0 (f|_\Vspace) ]^C ,
    \label{eqn:DGM-0} \\
  \Ddgm_p (f|_\Mspace)  &=&  \Ddgm_p (f|_\Uspace) \sqcup \Ddgm_p (f|_\Vspace) ,
    \label{eqn:DGM}   \\
  \Ddgm_n (f|_\Mspace)  &=&  [ \Ddgm_0 (f|_\Uspace) \sqcup \Ddgm_0 (f|_\Vspace) ]^{CT} ,
    \label{eqn:DGM-n}
\end{eqnarray}
for $1 \leq p \leq n-1$, 
where $C$ stands for applying the cascade
and $T$ stands for reflecting the diagram.
\end{result}
We note that the assumption of $n$ being positive is necessary
since the formulas do not hold for $n = 0$.
Concerning the proof of the theorem, 
note that \eqref{eqn:DGM} is clear, and  that \eqref{eqn:DGM-n} follows from
\eqref{eqn:DGM-0} using the duality of persistence diagrams from~\cite{CEH09}.
We will need to study the impact of the non-trivial groups $\Sgroup_i^0$
to prove \eqref{eqn:DGM-0}.

\paragraph{Latitudinal manifolds.}
Call the minimum of $f$ the \emph{south-pole}
and the maximum the \emph{north-pole} of the sphere.
Let $\Mspace'$ be a component of $\Mspace$, and note that it
is an $n$-manifold that decomposes $\Sspace$ into two complementary subsets.
If it separates the two poles, we refer to $\Mspace'$
as a \emph{latitudinal $n$-manifold}.
Assuming neither pole lies on $\Mspace$,
we order the latitudinal $n$-manifolds from south to north as
$\Mspace_1, \Mspace_2, \ldots, \Mspace_\ell$.
Letting $u_k$ and $w_k$ be the minimum and maximum values of $f$
restricted to $\Mspace_k$, we get
$u_1 < u_2 < \ldots < u_\ell$ as well as $w_1 < w_2 < \ldots < w_\ell$.
 
\begin{figure}[hbt]
  \vspace{0.0in}
  \centering
  \resizebox{!}{1.2in}{\input{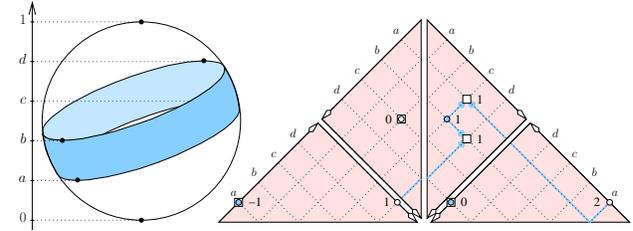}}
  \caption{Left: the height function on the $2$-sphere decomposed into
    an annulus and a pair of disks.
    Right: the reduced persistence diagram of $f$ restricted to
    the annulus (shaded circular dots), to the pair of disks
    (white circular dots), and to the two circles (squares).}
  \label{fig:counterexample}
\end{figure}

For each component $\Mspace'$ of $\Mspace$, there are
\emph{neighboring components} $\Uspace'$ of
$\Uspace$ and $\Vspace'$ of $\Vspace$
defined such that $\Uspace' \capsp \Vspace' = \Mspace'$.
A single component of $\Uspace$ or $\Vspace$
can be neighbor to an arbitrary number
of $n$-manifolds, but not to more than two latitudinal $n$-manifolds.
Specifically, there are components
$\Sspace_0, \Sspace_1, \ldots, \Sspace_\ell$ of $\Uspace$ and $\Vspace$
such that $\Mspace_k = \Sspace_{k-1} \capsp \Sspace_k$
for each $1 \leq k \leq \ell$.
We refer to the $\Sspace_k$ as \emph{latitudinal components}.
For example, in Figure \ref{fig:counterexample},
we have $\ell = 2$ latitudinal $1$-manifolds
and $\ell + 1 = 3$ latitudinal components.
Setting $u_0 = 0$ and $w_{\ell+1} = 1$,
the minimum and maximum values of $f$ restricted
to $\Sspace_k$ are $u_k$ and $w_{k+1}$, for $0 \leq k \leq \ell$.
Assuming $\Sspace_k$ belongs to $\Uspace$, it gives rise to zero or more
dots in the ordinary subdiagram of $\Ddgm_0 (f|_\Uspace)$ and to
exactly one dot, $(u_k, w_{k+1})$, in the horizontal subdiagram.
We say the dot in the horizontal subdiagram \emph{represents}
the $0$-dimensional homology class defined by $\Sspace_k$.
Note that its coordinates are indexed consistently with the
notation used in the introduction of the cascade.
There is indeed a connection, namely the dots in the
$0$-th horizontal subdiagram representing latitudinal
components are the extreme ones in the multiset:
\begin{result}[Extrema Lemma]
  The latitudinal components of $\Uspace$ and $\Vspace$
  correspond bijectively to the extreme
  dots in the disjoint union of
  $\Ddgm_0 (f|_\Uspace)$ and $\Ddgm_0 (f|_\Vspace)$.
\end{result}
We omit the proof, which is not difficult.
The statement includes the case in which there is no latitudinal $n$-manifold
so that $\Sspace_0$ is the only latitudinal component.
It contains both poles and is therefore represented by $(0,1)$,
which is the only extreme dot in the disjoint union of $0$-th diagrams.

\paragraph{Proof of \eqref{eqn:DGM-0}.}
We compare $\Ddgm_0(f|_\Mspace)$ with the disjoint union of
 $\Ddgm_0(f|_\Uspace)$ and $\Ddgm_0 (f|_\Vspace)$,
 noting that their dots all belong to the ordinary and horizontal
 subdiagrams.
 Consider first a dot $(a,b)$ in the ordinary subdiagram of $\Ddgm_0 (f|_\Mspace)$.
 It represents a component in the sublevel set that is born at
 a minimum $x \in \Mspace$, with $f(x) = a$, and that dies at a saddle point
 $y \in \Mspace$, with $f(y) = b$.
 Let $\Mspace' \subseteq \Mspace$ be the connected $n$-manifold that
 contains $x$ and $y$, and recall that $\Mspace' = \Uspace' \capsp \Vspace'$.
 Assume the neighboring component $\Uspace'$ of $\Uspace$ lies above $x$.
 In the sequence of sublevel sets of $f|_\Uspace$, we see the birth
 of a component at $f(x) = a$ and its death at $f(y) = b$.
 It follows that $(a,b)$ is also a dot in the ordinary subdiagram
 of $\Ddgm_0 (f|_\Uspace)$ and therefore of the disjoint union of the
 diagrams of $f|_\Uspace$ and $f|_\Vspace$.
 The argument can be reversed, which implies that the
 $0$-th ordinary subdiagrams are the same.

 Consider second a dot $(a,b)$ in the horizontal subdiagram of
 $\Ddgm_0 (f|_\Mspace)$.
 It represents a connected $n$-manifold $\Mspace' \subseteq \Mspace$,
 which splits $\Sspace$ into two subsets.
 If $\Mspace'$ is non-latitudinal, then one subset contains both poles
 and the other contains neither.
 The latter subset contains a (non-latitudinal) neighboring component,
 which is represented by $(a,b)$ in the horizontal subdiagram of
 $\Ddgm_0 (f|_\Uspace)$ or of $\Ddgm_0 (f|_\Vspace)$.
 Again, the argument can be reversed.
 If on the other hand, $\Mspace'$ is latitudinal, then
 $(a,b) = (u_i,w_i)$, for some $i$,
 where we write $u_1, u_2, \ldots, u_\ell$ for the minimum values
 and $w_1, w_2, \ldots, w_\ell$ for the maximum values of the
 latitudinal $n$-manifolds, as before.
 More generally, we get the dots $(u_i, w_i)$ in the horizontal subdiagram
 of $\Ddgm_0 (f|_\Mspace)$, for $1 \leq i \leq \ell$.
 Adding $u_0 = 0$ and $w_{\ell+1} = 1$, we get the dots
 $(u_i, w_{i+1})$ in the horizontal subdiagram of
 $\Ddgm_0 (f|_\Uspace) \sqcup \Ddgm_0 (f|_\Vspace)$, for $1 \leq i \leq \ell$.
 After adding $(0,0)$ and $(1,1)$, which are both irrelevant,
 we get precisely the dots specified in
 \eqref{eqn:cascade-new} and in \eqref{eqn:cascade-old}.
 This implies that the two diagrams are related to each other
 by a cascade, which completes the proof of \eqref{eqn:DGM-0}.

\begin{figure}[hbt]
 \vspace{0.0in}
 \centering
 \resizebox{!}{0.95in}{\input{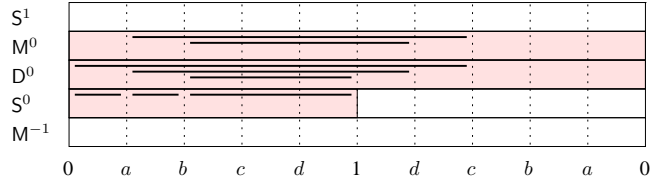}}
 \caption{White areas represent trivial groups.
   The pattern of having every third row white is broken
   at the top and at the bottom.
   The black life-time intervals show homology classes
   that appear in the analysis of the $0$-th and $n$-th
   persistence diagrams.}
 \label{fig:pattern}
\end{figure}

As an alternative to the above geometric proof of \eqref{eqn:DGM-0},
we could give an algebraic argument that considers the filtrations
of Meyer-Vietoris sequences; see Figure \ref{fig:pattern},
which draws the filtrations as rows in a matrix,
shading the area of non-trivial groups.
The obstacle to $\Mgroup_i^0$ and $\Dgroup_i^0$ being isomorphic,
and thus the reason for the cascade,
are the non-trivial groups $\Sgroup_i^0$ in the left half of the filtration.
Because each non-trivial $\Sgroup_i^0$ has exactly one generator, 
$\Dgroup_i^0$ must have one more generator than $\Mgroup_i^0$.
Indeed, by comparing the groups, before and after the cascade,
we see that in the left half,
each $\Dgroup_i^0$ has one more generator than $\Mgroup_i^0$.

\ignore{
\paragraph{Reduced case.}
With the General Shore Theorem, 
we have arrived at a complete description of the diagram of $f|_\Mspace$
in terms of the diagrams of $f|_\Uspace$ and $f|_\Vspace$.
We also state the relations for the reduced diagram, which are
\begin{eqnarray*}
  \rDdgm_p (f|_\Mspace)  &=&  \rDdgm_p (f|_\Uspace) \sqcup \rDdgm_{n-p} (f|_\Uspace)^T ,\\
  \rDdgm_n (f|_\Mspace)  &=&  [ \rDdgm_0 (f|_\Uspace) \sqcup \rDdgm_n (f|_\Uspace)^T ]^{CT} ,
\end{eqnarray*}
for $0\leq p < n$.
The result is different for two reasons.
First, all $0$-th reduced homology groups of sublevel sets of $f$  are trivial,
so there is no need for the cascade.
Second, we can use the Land and Water Theorem to
substitute the reflected reduced diagram of $f|_\Uspace$
for the reduced diagram of $f|_\Vspace$.
}

\paragraph{Euclidean case.}
Considering that \eqref{eqn:MVADSphere} to \eqref{eqn:MVADopenBall}
give elegant relations between the Betti numbers of the shore
and the land, it is perhaps surprising that we need cascades to
formulate a similar result for persistence diagrams.
Indeed, \eqref{eqn:MVADSphere} to \eqref{eqn:MVADopenBall} suggest
that the persistence diagrams of $f|_\Mspace$ be the disjoint union
of the persistence diagram of $f|_\Uspace$ and of its reflection.
The example in Figure \ref{fig:counterexample} shows that this
simple relation does not hold in general.
Indeed, we have the dots $(a,c)$, $(b,d)$, $(c,a)$, $(d,b)$
in the standard diagram of $f|_\Mspace$,
while the standard diagram of $f|_\Uspace$ contains $(a,d)$ and $(c,b)$.
This clearly violates the suggested relation.
However, there is a natural setting in which the relation does hold,
which we now describe.
\begin{result}[Euclidean Shore Theorem]
  Let $n$ be a positive integer,
  let $\Aspace$ be a compact set
  whose boundary is an $n$-manifold in $\Rspace^{n+1}$,
  and suppose that $e: \Rspace^{n+1} \to \Rspace$ has no
  homological critical values 
  and its restriction to $\partial \Aspace$ is tame.
  Then $\Ddgm (e|_{\partial \Aspace}) 
    = \Ddgm (e|_\Aspace) \sqcup \Ddgm (e|_\Aspace)^T$.
\end{result}
\proof
 We can extend $e$ to a perfect Morse function
 $f: \Sspace^{n+1} \to \Rspace$ and $\Aspace$ to a subset
 $\Uspace$ of $\Sspace = \Sspace^{n+1}$ such that the
 persistence diagrams of $e$ restricted to $\Aspace$
 and to $\partial \Aspace$ are the same as those of $f$ restricted
 to $\Uspace$ and $\Mspace = \partial \Uspace$.
 It thus suffices to show that the persistence diagram of $f|_\Mspace$
 is the disjoint union of the diagram of $f|_\Uspace$ and of its reflection.
 For $0<p<n$, this follows from \eqref{eqn:DGM}, from
 \begin{eqnarray*}
   \rDdgm_p(f|_\Vspace) &=& \rDdgm_{n-p}(f|_\Uspace)^T ,
 \end{eqnarray*}
 as stated in the Land and Water Theorem,
 and the fact that the reduced diagrams are equal to the non-reduced ones.
 For $p = 0$, we start with $\eqref{eqn:DGM-0}$ and note that there is only one
 extreme component in $\Ddgm_0 (f|_\Uspace) \sqcup \Ddgm_0 (f|_\Vspace)$,
 namely the one of $\Vspace$ that contains
 both the minimum and the maximum of $f$.
 It follows that the cascade leaves $\Ddgm_0 (f|_\Uspace)$
 unchanged while it turns $\Ddgm_0 (f|_\Vspace)$
 into $\rDdgm_0 (f|_\Vspace)$. Using the Land and Water Theorem again, we
 obtain:
 \begin{eqnarray*}
   \Ddgm_0 (f|_\Mspace) &=& [\Ddgm_0 (f|_\Uspace) \sqcup \Ddgm_0 (f|_\Vspace)]^{C}  \\
                     &=& \Ddgm_0 (f|_\Uspace) \sqcup \rDdgm_0 (f|_\Vspace)   \\
                     &=& \Ddgm_0 (f|_\Uspace) \sqcup \rDdgm_n (f|_\Uspace)^T \\                     &=& \Ddgm_0 (f|_\Uspace) \sqcup \Ddgm_n (f|_\Uspace)^T.
 \end{eqnarray*}
 Finally, for $p = n$, we exploit Poincar\'{e} duality for manifolds,
 which implies $\Ddgm_n (f|_\Mspace) = \Ddgm_0 (f|_\Mspace)^T$;
 see \cite{CEH09}.
\eop

The $4$-dimensional version of the theorem, for $n = 3$,
brings us back full circle to the motivation
for this work, namely the computation of the persistence diagram of
the space-time shape formed by a moving collection of biological cells \cite{EHKKS11}.
Modeling space-time as $\Rspace^4$,
the shape is compact, and we consider the time function
restricted to that shape.
This data satisfies the assumptions of the Euclidean Shore Theorem,
so we can infer the persistence diagram of the function on the boundary from
the diagram of the function on the solid $4$-dimensional shape.

\section{Discussion}
\label{sec6}

The main contributions of this paper
are two extensions of Alexander duality from spaces to functions.
The first extension is direct and relates the persistence diagrams
of a perfect Morse function restricted to two complementary subsets
of the $(n+1)$-dimensional sphere with each other.
The second extension relates
the persistence diagram of the function restricted
to the intersection of the complementary subsets
with the diagram of the function restricted to one subset.
A key tool in its proof is the filtration of Mayer-Vietoris sequences
(or the Mayer-Vietoris sequence of filtrations).
This suggests the study of more general filtrations of exact sequences.
Besides the hope to develop a general purpose device that pervades
persistent homology in the same way exact sequences pervade homology,
we motivate the study with a few concrete questions.
\begin{itemize}
  \item  Is it possible to relate the extent to which a function
    $f$ on $\Sspace^{n+1}$ is not perfect Morse with the severity
    with which this function violates our two extensions of Alexander duality?
  \item  Does every Morse function $f$ on $\Sspace = \Sspace^{n+1}$
    have a decomposition into complementary subsets, $\Sspace = \Uspace \cupsp \Vspace$,
    such that up to small modifications,
    the disjoint union of the diagrams of $f$ and $f|_\Mspace$
    is equal to the disjoint union of the diagrams of
    $f|_\Uspace$ and $f|_\Vspace$?
  \item  Can these or similar relations be developed into
    a divide-and-conquer algorithm for computing the
    persistence diagram of $f$ on a sphere?
\end{itemize}
It would furthermore be interesting to generalize the results
of this paper, as well as the above questions, to functions
on $(n+1)$-manifolds other than the sphere.


\end{document}